# Is a deterministic universe logically consistent with a probabilistic Quantum Theory?

*Bhupinder Singh Anand*

*If we assume the Thesis that any classical Turing machine T, which halts on every n-ary sequence of natural numbers as input, determines a PA-provable formula, whose standard interpretation is an n-ary arithmetical relation $f(x_1, ..., x_n)$ that holds if, and only if, T halts, then standard PA can model the state of a deterministic universe that is consistent with a probabilistic Quantum Theory. Another significant consequence of this Thesis is that every partial recursive function can be effectively defined as total.*

**1.0 Introduction**

In this paper, we consider some significant consequences of the two meta-theses:

(a) *Classical Halting Thesis*: Every PA[1]-unprovable formula, that is a true[2] n-ary arithmetical relation $f(x_1, ..., x_n)$ under the standard interpretation[3], determines a

---

[1] By "PA", we mean a classical, standard, first order formal system of Arithmetic such as Mendelson's formal system S (3, p102). We note, however, that the arguments of this paper would also hold for Gödel's formal system P (2, p9).

[2] We follow Tarski's definitions of the "satisfiability" and "truth" of a formula under a given interpretation (3, p49-53).

We note that Tarski's definitions implicitly assume that a formula [F] is true under an interpretation if, and only if, there is an effective method (native to the interpretation) for determining that the formula [F] is satisfied under the interpretation for all sequences in the domain of the interpretation (1, Appendix 5). It follows by Church's Thesis (3, footnote on p147), that, if the PA-formula [F] is true under the standard interpretation, then it is recursive.

[3] We follow Mendelson's definitions of interpretation (3, p49) and standard interpretation (3, 107).

classical Turing machine[4] T that halts on every $n$-ary sequence $k_1, ..., k_n$ of natural numbers as input if, and only if, $f(k_1, ..., k_n)$ holds.

(b) *Quantum Halting Thesis*: Any classical Turing machine T that halts on every $n$-ary sequence of natural numbers as input determines a PA-provable formula, whose standard interpretation is an $n$-ary arithmetical relation $f(x_1, ..., x_n)$ such that, for any sequence $k_1, ..., k_n$ of natural numbers, $f(k_1, ..., k_n)$ holds if, and only if, T halts.

Since the conclusion in (*a*), and the premise in (*b*), can never be effectively verified, these meta-theses are essentially undecidable; clearly, the two are mutually inconsistent.

We thus have two, apparently equally sound, but essentially different, systems of Arithmetic based on standard PA that should, reasonably, have significantly different consequences.

That this is, indeed, so, is seen in the following:

*Theorem 1*: A deterministic universe is inconsistent with the Classical Halting Thesis.

*Theorem 2*: A deterministic universe is logically consistent with the Quantum Halting Thesis and, ipso facto, with the probabilistic formulation of Quantum Theory.

*Theorem 3*: The Quantum Halting Thesis implies that every partial recursive function can be effectively defined as a total function.

---

[4] We follow Mendelson's description of a classical Turing machine (3, p229).

## 2.0 Defining a deterministic universe

*Definition 1*: The standard interpretation of a unary PA-formula is a quantum state if, and only if, the formula is PA-unprovable but true under the standard interpretation[5].

The significance of a quantum state $[F(x)]$[6] is that, given any natural number $k$, we can always construct an individually[7] effective[8] method for determining that $F(k)$ holds, but there is no uniformly effective method such that, given any $k$, we can always determine that $F(k)$ holds.

*Premise*: We assume, next, that, given a (physical) particle $P$, every property of the particle that can be expressed by any finite string in any recursively enumerable language[9] L, which has a finite alphabet, can be uniquely Gödel-numbered[10] by a

---

[5] In his seminal 1931 paper (2), Gödel constructs such a formula in an intuitionistically unobjectionable manner.

[6] We denote the PA-formula whose interpretation is $F$ by $[F]$.

[7] The thesis underlying the arguments in this paper is that we can differentiate between intuitive concepts such as "individually decidable", "individually computable", "individually terminating routine" etc., and their corresponding, and contrasting, intuitive concepts such as "uniformly decidable", "uniformly computable", "uniformly terminating routine" etc. (cf. 1, Appendices 4 & 5).

Loosely speaking, we would use the terminology, for instance "individually decidable", if, and only if, say, a PA-formula with free variables is decidable individually, as satisfiable or not, in an interpretation M for any given instantiation. This corresponds to, but is not necessarily a consequence of, the case where every instantiation of the formula is PA-provable. We note that the closure of the formula may not be P-provable

In contrast, we would use the terminology "uniformly decidable" if, and only if, a PA-formula with free variables is decidable jointly, as true or not, in the interpretation M as an infinitely compound assertion of all of its instantiations. This corresponds to, and is a consequence of, the case where the closure of the formula is PA-provable.

The concept of an intuitively terminating routine may be taken as corresponding to the classical concept of an "effective procedure" (3, p207).

[8] We broadly follow Mendelson's approach to "effectiveness" (3, p207).

[9] We define a language L as recursively enumerable if there exists a Turing machine that accepts every string of the language, and does not accept strings that are not in the language.

natural number. We assume, similarly, that every value of any property, which can be finitely expressed similarly in L, can also be uniquely Gödel-numbered by a natural number.

We note that not every natural number need correspond to a property, and that not every natural number need correspond to the value of a property. We also note that, even if the number of all possible particles, properties and values in a deterministic universe are not denumerable, we need only consider a denumerable number of particles, properties and values, since we can only express denumerable strings within L.

*Definition 2*: We define the state of a particle $P$ at any instant $t$[11] as a number-theoretic relation $p_t(x, y)$[12] that holds if, and only if, for any given property represented by $x$, there is a unique natural number $y$, where the domain of $x$ is assumed to be a recursively enumerable set[13] of the Gödel-numbers of all properties that are expressible in L, and the domain of $y$ is assumed to be a recursively enumerable set of the Gödel-numbers of all values that are expressible in L.

*Definition 3*: We define the universe as deterministic if, for any given particle $P$, any given property[14] $k$, and any given instant $t$, there is an effective method to determine a unique value $m$ such that $p_t(k, m)$ holds.

---

[10] By "Gödel-numbered", we mean the assignment of a unique natural number to the formal expression of any concept within a language (3, p135).

[11] We note that $t$ here is a natural number that is assumed to be the $t$'th measurement of the property of the particle $P$ with reference to some initial quantum state $p_0(x, y)$. In other words, we assume that, for any property $k$, the quantum states $p_0(k, m_0)$, $p_0(x, m_1)$, ..., $p_0(x, m_{t-1})$ are already known.

[12] We note that $p_t(x, y)$ may, in turn, be definable by a number-theoretic relation $P(x, y, z)$, where $z$ is the Gödel-number of any expression of the time $t$ in the language L.

[13] " A set of natural numbers is called recursively enumerable (r.e) if, and only if, it is either empty or is the range of a recursive function. Intuitively, if we accept Church's Thesis, then a recursively enumerable set is a collection of natural numbers which is generated by some mechanical process." (3, p250)

*Lemma 1*: If the properties expressible in L are recursively enumerable, and we assume Church's Thesis[15], then we can always effectively define $p_t(k, m)$ as holding if $k$ is not the Gödel-number of a property. Hence $p_t(x, y)$ is total.

*Proof*: By definition, and the Church Thesis, there is a classical Turing machine T such that, input any given natural numbers $k, m$, T halts if, and only if, $p_t(k, m)$ holds. If $k$ is not the Gödel-number of a property, then $p_t(k, m)$ is undefined, and we assume that, in such a case, T will loop on the input $k$.

We note that the definition of any classical Turing machine T can be extended to allow for recording, in its infinite memory[16], of every instantaneous tape description[17] during the operation of a classical Turing machine.

Now, T can always be meta-programmed to compare its current instantaneous tape description with the finite set of previous instantaneous tape descriptions during its current operation, and to self-terminate if an instantaneous tape description repeats itself, i.e. at the onset of a loop[18].

We can, thus, use this condition to effectively define $p_t(k, m)$ as holding if $k$ is not the Gödel-number of a property. Hence $p_t(x, y)$ is total.

---

[14] Where there is no obvious ambiguity, we shall henceforth use the terms "property $k$" and "value $m$" when we actually mean "property whose Gödel-number is $k$", and "value whose Gödel-number is $m$".

[15] We take Church's Thesis to essentially state that a number-theoretic function is effectively computable if, and only if, it is recursive, and a partial number-theoretic function is effectively computable if, and only if, it is partial recursive (3, p227). We note too, that a number-theoretic function is Turing-computable if, and only if, it is partial recursive (3, p233, Corollary 5.13 & p237, Corollary 5.15). We thus have the Turing Thesis (4) that a number-theoretic function is effectively computable if, and only if, it is Turing-computable.

[16] We may visualise T as a virtual teleprinter that copies, and records, every instantaneous tape description in an area of the machine that acts as a dynamic, random access, memory during any operation of T.

[17] Cf. Mendelson (3, p230).

[18] Cf. Anand (1, Appendix 7).

*Lemma 2*: If we assume the Church Thesis[19], it follows from this definition that, in a deterministic universe, the number-theoretic relation $(E!y)p_t(x, y)$[20], describing the complete state of the particle $P$ at instant $t$, is recursive[21].

*Proof*: This follows since every Turing-computable function $f$ is recursive if it is total[22].

We now consider the case where the PA-formula $[Q_t(x, y)]$, which expresses[23] the recursive relation $p_t(x, y)$ in PA, is such that the formula $[(E!y)Q_t(x, y)]$ is a quantum state, i.e. it is PA-unprovable, but true under the standard interpretation of PA. It follows that, given any natural numbers $k, m$:

(*i*)  $(E!y)Q_t(k, y)$ holds, and

(*ii*)  $Q_t(k, m)$ holds if, and only if, $p_t(k, m)$ holds.

### 3.0 Determinism and the Classical Halting Thesis

We now have that:

*Lemma 3*: If we assume the Classical Halting Thesis (*a*), then there is some classical Turing machine T such that, for any input $k$, T will halt and return the value $m$.

---

[19] We can also assume, alternatively, both the Restricted Church Thesis and the Universal Church Thesis (1, Appendix 5); prima facie, the Quantum Halting Thesis appears to be a trivial consequence of these. However, this may not, then, highlight the contrast between the Classical and Turing Halting Theses as is intended in this paper.

[20] The notation "(E!$x$)" denotes uniqueness.

[21] We follow Mendelson's (3, p120) definitions of recursive functions and relations.

[22] Cf. Mendelson (3, p233, Corollary 5.13).

[23] We follow Mendelson's (3, p117-118) definitions of expressibility and representability of number-theoretic functions.

*Proof*: T will halt on input $k$ if, and only if, $(E!y)Q_t(k, y)$ holds. By our hypothesis, for any given $k$, $(E!y)Q_t(k, y)$ holds uniquely for some value $m$. Hence, for given input $k$, T must run a sub-routine that will halt on, and can be programmed to return, $m$.

This would clearly imply that:

*Lemma 4*: If we assume the Classical Halting Thesis (*a*), then there is some experiment that can completely determine the state of any finite number of, arbitrarily selected, properties simultaneously. However, this would violate Quantum Uncertainty, which postulates states of a particle that cannot be completely determined simultaneously.

We thus conclude:

*Theorem 1*: A deterministic universe is inconsistent with the Classical Halting Thesis.

**4.0 Determinism and the Quantum Halting Thesis**

However, it is obvious that:

*Lemma 5*: If we assume the Quantum Halting Thesis (*b*), then there is no classical Turing machine T such that, for any given input $k$, T will halt and return the value $m$, since this would make the formula $[(E!y)Q_t(x, y)]$ PA-provable, contrary to our hypothesis.

This, now, clearly implies that:

*Lemma 6*: There is no experiment, or finite set of experiments, that can completely determine the state of a finite number of, arbitrarily selected, properties of a particle $P$ simultaneously, even though, by our premise, given any property $k$, there is always some experiment that will completely determine its state $m$ at instant $t$.

We thus conclude that:

*Theorem 2*: A deterministic universe is logically consistent with the Quantum Halting Thesis, and ipso facto, with the probabilistic formulation of Quantum Theory[24].

## 5.0 Effectively defining every partial recursive function as total

Curiously, the above consequences emerged incidentally out of the following argument, which effectively defines every partial recursive function as total.

*Theorem 3*: The Quantum Halting Thesis implies that every partial recursive[25] number-theoretic function $f(x_1, ..., x_n)$ is constructively definable as a total function.

*Proof*: We assume that $f$ is obtained from the recursive function $g$ by means of the unrestricted $\mu$-operator, so that $f(x_1, ..., x_n) = \mu y(g(x_1, ..., x_n, y) = 0)$.

Given any sequence $k_1, ..., k_n$ of natural numbers, we can then define:

(*i*) a classical Turing machine T1 that will halt if, and only if, $(g(k_1, ..., k_n, k) = 0)$ holds for some natural number k;

---

[24] We note that, on this model, the probabilities of Quantum Theory can be interpreted as describing, not any essential property pertaining to the physical state of a particle, but the probability that any relation $p_t(x, y)$ that is satisfied by any experiment, or finite series of experiments, describes the complete physical state of the particle at the instant $t$.

[25] Classically (3, p120, 121, 214), a partial function $f$ of $n$ arguments is called partial recursive if, and only if, $f$ can be obtained from the initial functions (zero function), projection functions, and successor function (of classical recursive function theory) by means of substitution, recursion and the classical, unrestricted, $\mu$-operator. $f$ is said to come from $g$ by means of the unrestricted $\mu$-operator, where $g(x_1, ..., x_n, y)$ is recursive, if, and only if, $f(x_1, ..., x_n) = \mu y(g(x_1, ..., x_n, y) = 0)$, where $\mu y(g(x_1, ..., x_n, y) = 0)$ is the least number $k$ (if such exists) such that, if $0 =< i =< k$, $g(x_1, ..., x_n, i)$ exists and is not 0, and $g(x_1, ..., x_n, k) = 0$. We note that, classically, $f$ may not be defined for certain $n$-tuples; in particular, for those $n$-tuples $(x_1, ..., x_n)$ for which there is no $y$ such that $g(x_1, ..., x_n, y) = 0$.

(ii) a classical Turing machine T2 that will halt if, and only if, $\sim(g(k_1, ..., k_n, k) = 0)$ holds for some natural number k;

(iii) a classical Turing machine T3[26] that will halt if, and only if, $[H(k_1, ..., k_n, y)]$ is PA-provable, where $[H(x_1, ..., x_n, y)]$ is the PA-formula that represents $(g(x_1, ..., x_n, y) = 0)$ in PA.

If we now assume the Quantum Halting Thesis, it follows that:

(a) if T1 halts, without looping, on all natural numbers, then $[H(k_1, ..., k_n, y)]$ is PA-provable;

(b) similarly, if T2 halts, without looping, on all natural numbers, then $[H(k_1, ..., k_n, y)]$ is PA-provable.

Hence, if $(g(k_1, ..., k_n, y) = 0)$ holds for all natural numbers y, then:

(c) either $[H(k_1, ..., k_n, y)]$ is PA-provable, and T3 will halt;

(d) or, if $[H(k_1, ..., k_n, y)]$ is not PA-provable, then T2 will loop, and self-terminate, for some natural number $k$.

If we now run T1, T2 and T3 simultaneously, then we are assured that at least one of them will halt, or self-terminate, for some natural number $k$. We can thus use the condition indicated by the halting, or self-terminating, state to effectively define the partial recursive function $f(x_1, ..., x_n)$ appropriately as a total partial recursive function.

---

[26] The definition of T3 is based on Gödel's recursive number-theoretic relation $xBy$ (2, Definition 45), which holds if, and only if, $x$ is the Gödel-number of a PA-proof of the PA-formula whose Gödel-number is $y$.

**Appendix 1: Implications of Definition 2**

We note that *t* in *Definition 2* is a natural number that is assumed to be the *t*'th measurement of the property of the particle *P* with reference to some initial quantum state $p_0(k, m_0)$. In other words, we assume that, for any property $k$, the quantum states $p_0(k, m_0), p_0(k, m_1), ..., p_0(k, m_{t-1})$ are already known.

We note that, for any given property $k$ and instant $t = t_n$, we can define a recursive quantum Gödel *ß*-function[27] $q_{P, k}(b, c, i)$, such that $q_{P, k}(b, c, i) = m_i$ for all $0 =<$[28] $i =< n$, where $m_i$ is the value of the property $k$ at the instant $t_i$, and $b$, $c$ are natural numbers that can be constructively determined by the sequence $m_0, m_1, ..., m_n$.

*Properties with random values in a deterministic universe*

If we assume that the instants $t_0, t_1, ..., t_n$ refer to a series of measurements that are taken in an experiment, then we may reasonably define the property $k$ as deterministic if the values represented by the sequence $m_0, m_1, ..., m_n$, as we increase the number of measurements $n$ indefinitely, are always assumed to be a sub-sequence of a Cauchy sequence; otherwise, we may, appropriately, define the property $k$ as random.

Prima facie, it thus follows that:

> *Lemma 7*: Both deterministic and random properties are consistent with a deterministic universe.

---

[27] We follow Mendelson's exposition (3, p131) of the Gödel *ß*-function.

[28] A notation such as "=<" is to be read naturally as "equal to or less than".

*Probability of a measurement yielding a given value*

We note, further, that, for any given sequence $m_0, m_1, ..., m_n$, there are denumerable pairs of natural numbers $(b, c)$ that define Gödel $\beta$-functions relative to the sequence; however, although they yield identical values for $0 =< i =< n$, each of these will yield different terms $m_{n+1}, m_{n+2}, ...$ for $i > n$. The challenge of any theory, thus, is to determine, firstly, which of these functions are such that all the terms of the non-terminating sequence are Gödel-numbers of values of the property $k$; secondly, the probability that the result $m_{n+1}$ of the measurement $t_{n+1}$, for any pair $(b, c)$, will represent any given value of the property $k$.

*Superposition*

We also note that, if $q_{P,k}(b, c, i) = m_i$ for all $0 =< i =< n$, where $m_i$ is the value of the property $k$ at the instant $t_i$, and $q'_{P,k}(b', c', i) = m'_i$ for all $0 =< i =< n'$, where $m'_i$ is the value of the property $k$ at the instant $t'_i$, we can always combine the two sequences of measurements $m_0, m_1, ..., m_n$, and $m'_0, m'_1, ..., m'_{n'}$, appropriately[29] to yield a sequence $m''_0, m''_1, ..., m''_{n''}$, where $n'' = n + n'$.

We thus have a Gödel $\beta$-function $q''_{P,k}(b'', c'', i) = m''_i$ for all $0 =< i =< n''$, where $m''_i$ is the value of the property $k$ at the instant $t''_i$, which may be considered as the superposition of the Gödel $\beta$-function $q'_{P,k}(b', c', i)$ on $q_{P,k}(b, c, i)$.

It follows that:

> *Lemma 8*: The probability[30] that $q''_{P,k}(b'', c'', n''+1) = m$ for a given $m$ is always equal to, or greater than, the probability that $q_{P,k}(b, c, n+1) = m$.

---

[29] We assume that there is an effective way of ensuring that the combined sequence $t''_i$, $0 =< i =< n''$, is also chronological.

*Proof*: The sequence $m_0, m_1, ..., m_n$ is a proper sub-sequence of $m''_0, m''_1, ..., m''_{n''}$. Hence, if $q''_{P, k}(b'', c'', n''+1) = m$, then there is always some pair $(b, c)$ such that $q_{P, k}(b, c, n''+1) = m$. The converse is not true.

Hence the set of Gödel *ß*-functions for which $q''_{P, k}(b'', c'', i) = m''_i$ for all $0 =< i =< n''$, where $m''_i$ is the value of the property $k$ at the instant $t''_i$, is a proper sub-set of the set of Gödel *ß*-functions for which $q_{P, k}(b, c, i) = m_i$ for all $0 =< i =< n$, where $m_i$ is the value of the property $k$ at the instant $t_i$.

*(Updated: Thursdy 2<sup>nd</sup> January 2003 4:10:45 AM by re@alixcomsi.com)*

---

[30] We assume here that this probability is, in a sense, "proportional" to the cardinality of the concerned set of Gödel *ß*-functions.